%% file: Partial-k-trees-arXiv.tex
\documentclass[runningheads,11points]{llncs}

\usepackage{amssymb,amstext}
\usepackage{cite,scrtime}
\usepackage[bookmarks=false]{hyperref}
\usepackage{graphics}
\usepackage{amsmath,dsfont}
\usepackage{mathrsfs}
\usepackage{fancyhdr}
\usepackage{verbatim}
\usepackage{boxedminipage}
\usepackage{colortbl}
\usepackage{tikz}

\usepackage[T1]{fontenc}
\usepackage{tikz}
\usepackage{graphicx}

\usepackage[lined,boxed,commentsnumbered]{algorithm2e}

\input{newcommand}

\tikzstyle{vertex}=[
  circle,
  draw,
  thick,
  inner sep=0.1cm,
  minimum size=2mm
]

\title{The number of labeled graphs\\of bounded treewidth
}

\author{Julien Baste~\inst{1}  \and Marc Noy~\inst{2} \and Ignasi Sau~\inst{1}}

\authorrunning{Julien Baste, Marc Noy, and Ignasi Sau}
\titlerunning{The number of labeled graphs of bounded treewidth}

\institute{CNRS, LIRMM, Universit\'e de Montpellier, Montpellier, France.\\
\email{\{baste,sau\}@lirmm.fr} \and
Department of Mathematics of Universitat Polit\`ecnica de Catalunya and\\Barcelona Graduate School in Mathematics, Barcelona, Catalonia.\\
\email{marc.noy@upc.edu}}

\begin{document}
\maketitle

\begin{abstract}
We focus on counting the number of labeled graphs on $n$ vertices and treewidth at most $k$ (or equivalently, the number of labeled partial $k$-trees), which we denote by $T_{n,k}$. So far, only the particular cases $T_{n,1}$ and $T_{n,2}$ had been studied. We show that
$$\left(c \cdot \frac{k\cdot 2^k \cdot n}{\log k} \right)^n \cdot  2^{-\frac{k(k+3)}{2}} \cdot k^{-2k-2}\ \leq\ T_{n,k}\ \leq\ \left(k \cdot 2^k \cdot n\right)^n \cdot 2^{-\frac{k(k+1)}{2}} \cdot k^{-k},$$
for $k > 1$ and  some explicit absolute constant $c > 0$. The upper bound is an immediate consequence of the well-known number of labeled $k$-trees, while the lower bound is obtained from an explicit algorithmic construction. It follows from this construction that both bounds also apply to graphs of pathwidth and proper-pathwidth at most $k$.

\end{abstract}

\section{Introduction}
\label{sec:intro}

Given an integer $k >0$, a \emph{$k$-tree} is a graph that can be constructed starting from a $(k+1)$-clique and iteratively adding a vertex connected to $k$ vertices that form a clique.
These graphs are natural extensions of trees, which correspond to $1$-trees, and received considerable attention since the late 1960s~\cite{BP69,M69,F71,HaPa68}. The notion of $k$-tree was introduced by Harary and Palmer~\cite{HaPa68}, and the number of labeled $k$-trees on $n$ vertices was first found by Beineke and Pippert~\cite{BP69}; cf. Moon~\cite{M69} and Foata~\cite{F71} for alternative proofs.

\begin{theorem}[Beineke and Pippert \cite{BP69}]
\label{thm:ktree}
There are
${n \choose k} (kn -k^2+1)^{n-k-2}$ many  $n$-vertex labeled $k$-trees.
\end{theorem}


A \emph{partial $k$-tree} is a subgraph of a $k$-tree. For two integers $n,k$ with $0 < k \leq n$, let $T_{n,k}$ denote the number of $n$-vertex labeled partial $k$-trees. While the number of $n$-vertex labeled $k$-trees is given by Theorem~\ref{thm:ktree}, it appears that very little is known about $T_{n,k}$.
Indeed, to the best of our knowledge, only the cases $k=1$ and $k=2$ have been  studied.
For the case $k=1$, the number of $n$-vertex labeled forests is asymptotically $T_{n,1} \sim \sqrt{e} \cdot n^{n-2}$~\cite{Forest90}.
For the case $k=2$, the number of $n$-vertex labeled series-parallel graphs, which is known to be exactly $T_{n,2}$, is asymptotically {$T_{n,2} \sim g \cdot n^{-\frac{5}{2}}\gamma^n n!$} for some explicit constants $g$ and $\gamma$~\cite{BGKN05}.

Interestingly, partial $k$-trees are exactly the graphs of \emph{treewidth} at most $k$. Treewidth is a structural graph invariant, which we formally define below, first introduced by Halin~\cite{Hal76} and later rediscovered by Robertson and Seymour~\cite{RoSe84} as a fundamental tool in their Graph Minors project culminating in the proof of Wagner's conjecture~\cite{RoSe04}.

A \emph{tree-decomposition} of width $k$ of a graph $G=(V,E)$ is a pair $({\sf T}, \mathcal{B})$, where ${\sf T}$ is a tree and $\mathcal{B} = \{ B_t \mid B_t \subseteq V, t \in V({\sf T}) \}$ such that:

\begin{enumerate}
\item $\bigcup_{t \in V({\sf T})} B_t = V$.
\item For every edge $\{u,v\} \in E$ there is a $t \in V({\sf T})$ such that $\{u, v\} \subseteq B_t$.
\item $B_i \cap B_{\ell} \subseteq B_j$ for all $\{i,j,\ell\} \subseteq V({\sf T})$ such that $j$ lies on the unique path from $i$ to $\ell$ in ${\sf T}$.
\item $\max_{t \in V({\sf T})} |B_t| = k +1$.
\end{enumerate}

The sets of $\mathcal{B}$ are called \emph{bags}. The \emph{treewidth} of $G$, denoted by $\tw(G)$, is the smallest integer $k$ such that there exists a tree-decomposition of $G$ of width $k$.
If ${\sf T}$ is a path, then $({\sf T}, \mathcal{B})$ is also called a \emph{path-decomposition}.
The \emph{pathwidth} of $G$, denoted by $\pw(G)$, is the smallest integer $k$ such that there exists a path-decomposition of $G$ of width $k$.

The following lemma is well-known and can be found, for instance, in~\cite{Kloks94}.
\begin{lemma}
\label{lemma:tw}
A graph has treewidth at most $k$ if and only if it is a partial $k$-tree.
\end{lemma}

Even if treewidth was introduced with purely graph-theoretic motivations, it turned out to have a number of algorithmic applications as well. One of the most relevant results in this area is
Courcelle's theorem~\cite{Courcelle90}, stating that any graph problem expressible in monadic second-order logic can be solved in linear time on graphs of bounded treewidth. Nowadays, treewidth is exhaustively used in both structural and algorithmic Graph Theory, cf. for instance the textbooks~\cite{Die10,CyganFKLMPPS15,Kloks94}. Recently, the treewidth of random graphs has also been studied under several probabilistic models~\cite{BoKl92, Gao12, LLO12, PeSe14, MiPe12}.



In this article, for any two integers $n,k$ with $0 < k \leq n$,
we are interested in counting the number of $n$-vertex labeled graphs that have treewidth at most $k$. By Lemma~\ref{lemma:tw}, this number is equal to $T_{n,k}$, and actually our approach relies heavily on the definition of partial $k$-trees.
As, by definition, the number of edges of an $n$-vertex $k$-tree is $kn - \frac{k(k+1)}{2}$,
by using Theorem~\ref{thm:ktree}
we obtain the following upper bound: 

\begin{equation}
\label{eq:ub}
T_{n,k}\ \leq\ 2^{kn - \frac{k(k+1)}{2}} \cdot {n \choose k} \cdot (kn -k^2+1)^{n-k-2}.
\end{equation}

\noindent Using the fact that ${n \choose k} \leq n^k$ and $1 \leq k^2$, from Equation~(\ref{eq:ub}) it follows that

\begin{equation}
\label{eq:sub}
T_{n,k}\ \leq\ (k \cdot 2^{k} \cdot n)^n \cdot 2^{ - \frac{k(k+1)}{2}} \cdot k^{-k}.
\end{equation}

On the other hand, we can easily provide a lower bound on $T_{n,k}$ with the following construction. Starting from an $(n-k+1)$-vertex forest, we add $k-1$ \emph{apices}, that is, $k-1$ vertices with an arbitrary neighborhood in the forest.
Every graph created in this way has exactly $n$ vertices and is clearly of treewidth at most $k$.
Moreover, as the number of labeled  forests on $n-k+1$ vertices is at least the number of trees on $n-k+1$ vertices, which is well-known to be $(n-k+1)^{n-k-1}$~\cite{Cay89}, and each apex can be connected to the forest in $2^{n-k+1}$ different ways, we obtain that
\begin{equation}
\label{eq:lb}
T_{n,k}\ \geq\ (n-k+1)^{n-k-1}\cdot 2^{(k-1)(n-k+1)}.
\end{equation}

As $n^{-2} \geq 2^{-n}$,
if we further assume that $\frac{n}{k}$ tends to infinity, from Equation~(\ref{eq:lb}) we get that, asymptotically,
\begin{equation}
\label{eq:slb}
T_{n,k}\ \geq\ \left(\frac{1}{4} \cdot 2^k\cdot n\right)^{n} \cdot 2^{-k^2}.
\end{equation}

The dominant factors of Equations~(\ref{eq:sub}) and~(\ref{eq:slb}), that is, $(k \cdot 2^{k} \cdot n)^n$ and $\left(\frac{1}{4} \cdot 2^k \cdot n \right)^n$, respectively, differ by a term $\left(\frac{1}{4}\right)^n$ and, most importantly, by a term $k^n$.


\bigskip

In order to close the gap between the existing lower and upper bounds on $T_{n,k}$, in this article we focus on improving the trivial lower bound presented above. We obtain the following result.

\begin{theorem}\label{thm:main}
For any two integers $n,k$ with $1 < k \leq n$, the number $T_{n,k}$ of $n$-vertex labeled graphs with treewidth at most $k$ satisfies
\begin{equation}
\label{eq:ilbi}
T_{n,k}\ \geq\ \left(\frac{1}{128e} \cdot \frac{k \cdot 2^{k} \cdot n }{ \log k}\right)^n\cdot 2^{-\frac{k(k+3)}{2}} \cdot k^{-2k -2}.
\end{equation}
\end{theorem}

That is, we fall short by a factor $(128e \cdot \log k)^n$ in order to reach the dominant factor of Equation~(\ref{eq:sub}).  In order to prove Theorem~\ref{thm:main}, we present in Section~\ref{sec:const} an algorithmic construction of a family of $n$-vertex labeled partial $k$-trees, which is inspired from the definition of $k$-trees. When exhibiting such a construction toward a lower bound, one has to play with the trade-off of, on the one hand, constructing as many graphs as possible and, on the other hand, being able to bound the number of duplicates; we perform this analysis in Section~\ref{sec:count}. Namely, we first count in Subsection~\ref{sub:numb} the number of elements created by the construction, and then we bound in Subsection~\ref{sub:bound} the number of times that  the same element may have been created. Finally, we present in Section~\ref{sec:further} some concluding remarks and several avenues for further research.

\section{The construction}
\label{sec:const}
Let $n$ and $k$ be two fixed positive  integers with $0 < k \leq n$.
In this section, we  construct a set $\Rcal_{n,k}$ of $n$-vertex labeled partial $k$-trees.
For notational simplicity, we let $R_{n,k} = |\Rcal_{n,k}|$. In Subsection~\ref{sub:def} we provide some notation and definitions used in the construction, in Subsection~\ref{sub:const}
we describe the construction, and in Subsection~\ref{sub:bounded}
 we prove that the treewidth of the generated graphs is indeed at most $k$. In fact,  we prove a stronger property, namely that the graphs we construct have \emph{proper-pathwidth} at most $k$, where the proper-pathwidth is a graph invariant that is lower-bounded by the pathwidth, which in turn is lower-bounded by the treewidth.


\subsection{Notation and definitions}
\label{sub:def}





For the construction, we use a \emph{labeling function} $\sigma$ defined by a permutation of $\dotsss{n}$ with the constraint that $\sigma(1) = 1$.
Inspired by the definition of $k$-trees, we will introduce vertices $\{v_1,v_2, \ldots, v_n\}$ one by one following the order $\sigma(1), \sigma(2), \ldots, \sigma(n)$ given by  the labeling function $\sigma$.
If $i,j \in \dotsss{n}$, then $i$ is called the \emph{index} of $\vv{i}$, the vertex $\vv{i}$ is the $i$-th introduced vertex
and,  if $j < i$, the vertex $\vv{j}$ is said to be \emph{to the left} of $\vv{i}$.

In order to build explicitly a class of partial $k$-trees, for every $i \geq k+1$ we will define:

\begin{enumerate}
\item
A set $A_i \subseteq \{j\mid j < i\}$ of  \emph{active} vertices, corresponding to the clique to which a new vertex can be connected in the definition of $k$-trees, such that $|A_i| = k$.
\item A vertex $a_i \in A_i$, called the \emph{anchor}, whose role will be described in the next paragraph.
\item An element $f(i) \in A_i$, called the \emph{frozen} vertex, which corresponds to a vertex that will not be active anymore.
\item A set $N(i) \subseteq A_i$, which corresponds to the indices of the neighbors of $\vv{i}$ to the left.
\end{enumerate}

The construction will work with \emph{blocks} of size $s$ for some integer $s$ depending of $n$ and $k$, to be specified in Subsection~\ref{sub:s}.  Namely, we will insert the vertices by consecutive blocks of size $s$, with the property that all vertices of the same block share the same anchor and are adjacent to it.

In the description of the construction,
we use the term \emph{choose} for the  elements
for which there are several choices, which will allow us to provide a lower bound on the number of elements in $\Rcal_{n,k}$.
It will be the case of the functions $\sigma$, $f$, and $N$.
As it will become clear later (cf. Section~\ref{sec:count}),
once $\sigma$, $f$, and $N$ are fixed, all the other elements of the construction are uniquely defined.

For every index $i$, we will impose that
$$|N(i)| > \frac{k+1}{2},$$
in order to have simultaneously enough choices for $N(i)$ and enough choices for the frozen vertex $f(i)$, which will be chosen among the vertices in $N(i-1)$.  On the other hand, as it will become clear in Subsection~\ref{sub:bound}, the role of the anchor vertices will be to uniquely determine the vertices belonging to ``its'' block.  For that, as we will see in the description of the construction, when a new block starts, its anchor is defined as the smallest currently active vertex.

\subsection{Description of the construction}
\label{sub:const}

Inspired by  the definition of  $k$-trees, we construct our partial $k$-trees in an algorithmic way. 
We say that a triple $(\sigma, f, N)$, with
$\sigma$ a permutation of $\dotsss{n}$,
$f : \{k+2, \ldots, n\} \rightarrow \dotsss{n}$, and
$N : \{2, \ldots, n\} \rightarrow 2^{\dotsss{n}}$,
is \emph{constructible} if it can be defined according to the following algorithm:

\medskip

\begin{algorithm}[H]
\IncMargin{1em}
\SetAlgoVlined
Choose $\sigma$, a permutation of $\dotsss{n}$ such that $\sigma(1) = 1$.

\For{i=2 \emph{\KwTo} k}{
 Choose $N(i) \subseteq \{{j}\mid j < i\}$, such that ${1} \in N(i)$.
}
\For{i=k+1}{
 Define $A_{k+1} = \{{j}\mid j < k+1\}$.\\
 Define $a_{k+1} = 1$.\\
 Choose $N(k+1) \subseteq \{{j}\mid j < i\}$, such that ${1} \in N(k+1)$.
}
\For{i=k+2 \emph{\KwTo} n}
{

\eIf{$i \equiv {k+2 \pmod s}$}{
 Define $f(i) = a_{i-1}$. \\
 Define $A_i = (A_{i-1} \sm \{f(i)\}) \cup \{i-1\}$.\\
 Define $a_i = \min A_i$.\\
 Choose $N(i) \subseteq A_{i}$ such that $a_i \in N(i)$ and $|N(i)| > \frac{k+1}{2}$; cf. Fig.~\ref{fig:introeq}.

}{
 Choose $f(i) \in (A_{i-1} \sm \{a_{i-1}\}) \cap N(i-1)$. \\
 Define $A_i = (A_{i-1} \sm \{f(i)\}) \cup \{i-1\}$.\\
 Define $a_i = a_{i-1}$.\\
 Choose $N(i) \subseteq A_{i}$ such that $a_i \in N(i)$ and $|N(i)| > \frac{k+1}{2}$; cf. Fig.~\ref{fig:intronoteq}.

}

}

\end{algorithm}





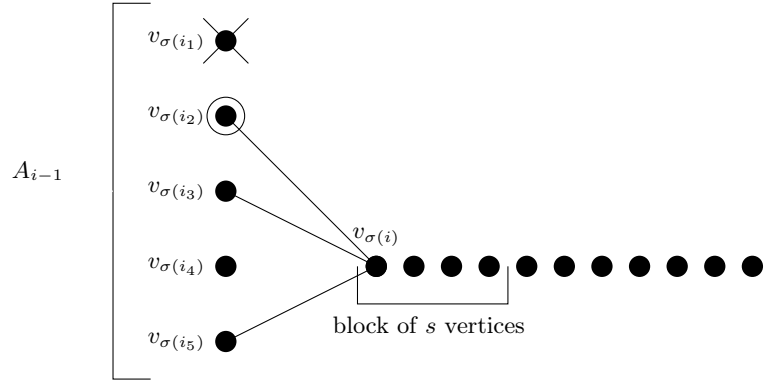
\begin{figure}[h!]
\centering
\begin{tikzpicture}
  \foreach \i in {1,...,5}
  {
    \node[place,label=left:$v_{\sigma(i_{\i})}$] at (0,5-\i) {};
  }

\draw (-1,4.5) -- (-1.5,4.5) -- (-1.5,-0.5) -- (-1,-0.5);
\capt{-3.5,2}{$A_{i-1}$};

    \node[place,label=above:$v_{\sigma(i)}$] at (2+0,1) {};
  \foreach \i in {0,...,10}
  {
    \node[place] at (2+0.5*\i,1) {};
  }
\draw (1.75,1) -- (1.75,.5) -- (3.75,.5) -- (3.75,1);

\forg{0,4}

\anch{0,3}

\draw (2,1) -- (0,3);
\draw (2,1) -- (0,2);
\draw (2,1) -- (0,0);

\capt{1.7,0}{block of $s$ vertices};
\end{tikzpicture}
  \caption{\ Introduction of $v_{\sigma(i)}$ with $k+2 \leq i \leq n$ and $i \equiv {k+2 \pmod s}$, $s = 4$, and $k = 5$.
    We assume that $i_1 < i_2 < i_3 < i_4 < i_5 < i$. We have defined $f(i) = \vv{i_1}$ and $a_i = \vv{i_2}$.
The frozen vertex $f(i)$ is marked with a cross, and the anchor $a_i$ is marked with a circle.
We choose $N(i) = \{i_3,i_5\}$.
}
  \label{fig:introeq}

\end{figure}

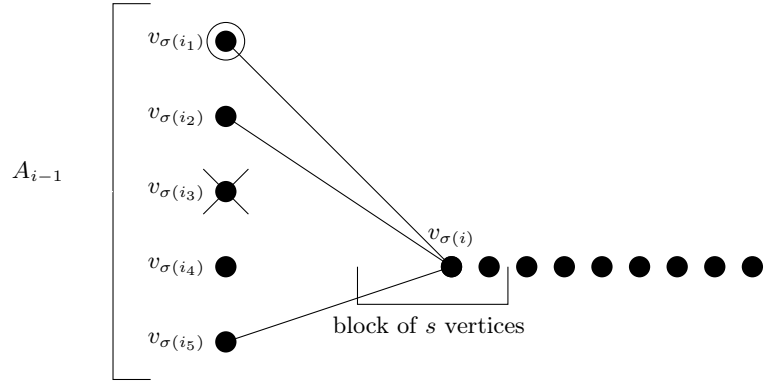
\begin{figure}[h!]
\centering
\begin{tikzpicture}
  \foreach \i in {1,...,5}
  {
    \node[place,label=left:$v_{\sigma(i_{\i})}$] at (0,5-\i) {};
  }
\draw (-1,4.5) -- (-1.5,4.5) -- (-1.5,-0.5) -- (-1,-0.5);
\capt{-3.5,2}{$A_{i-1}$};

    \node[place,label=above:$v_{\sigma(i)}$] at (3+0,1) {};
  \foreach \i in {2,...,10}
  {
    \node[place] at (2+0.5*\i,1) {};
  }
\draw (1.75,1) -- (1.75,0.5) -- (3.75,.5) -- (3.75,1);

\anch{0,4}

\forg{0,2}
\draw (3,1) -- (0,4);
\draw (3,1) -- (0,3);
\draw (3,1) -- (0,0);

\capt{1.7,0}{block of $s$ vertices};
\end{tikzpicture}
  \caption{\ Introduction of $v_{\sigma(i)}$ with $k+2 \leq i \leq n$ and $i \not \equiv {k+2 \pmod s}$, $s = 4$, and $k = 5$.
    We assume that $i_1 < i_2 < i_3 < i_4 < i_5 < i$. We have defined $a_i = a_{i-1} = \vv{i_1}$.
The frozen vertex $f(i)$ is marked with a cross, and the anchor $a_i$ is marked with a circle.
We choose  $f(i) = \vv{i_3}$, assuming $\vv{i_3}$ is a neighbor of $\vv{i_5}$, and  $N(i) = \{i_2,i_5\}$.
}
\label{fig:intronoteq}
\end{figure}

\medskip

 Let $(\sigma, f, N)$ be a constructible triple.
We define the graph $G(\sigma,f,N) = (V,E)$ such that
$V = \{v_i \mid i \in \dotsss{n}\}$ and $E = \{\{\vv{i}, \vv{j}\}\mid j \in N(i) \}$.
Note that, given  $(\sigma, f, N)$, the graph $G(\sigma,f,N)$ is well-defined.
We denote by $\Rcal_{n,k}$ the set of all graphs $G(\sigma, f, N)$ such that $(\sigma, f, N)$ is constructible.

\subsection{Bounding the proper-pathwidth of the constructed graphs}
\label{sub:bounded}

We start by defining the notion of proper-pathwidth of a graph. This parameter was introduced by Takahashi \emph{et al}.~\cite{TUK94}, and its relation with search games has been studied in~\cite{TUK95}.

Let $G$ be a graph and let $\mathcal{X} = \{X_1,X_2, \ldots, X_r\}$ be a sequence of subsets of $V(G)$.
The \emph{width} of $\mathcal{X}$ is $\max_{1 \leq i \leq r}|X_i|-1$. $\mathcal{X}$ is called a \emph{proper-path decomposition} of $G$ if the following conditions are satisfied:
\begin{enumerate}
\item For any distinct $i$ and $j$, $X_i \not \subseteq X_j$.
\item $\bigcup_{i=1}^r X_i = V(G)$.
\item For every edge $\{u,v\} \in E(G)$, there exists an $i$ such that $u,v \in X_i$.
\item For all $a$, $b$, and $c$ with $1 \leq a \leq b \leq c\leq r$, $X_a \cap X_c \subseteq X_b$.
\item For all $a$, $b$, and $c$ with $1 \leq a < b < c\leq r$, $|X_a \cap X_c| \leq  |X_b|-2$.
\end{enumerate}

The \emph{proper-pathwidth} of $G$, denoted by $\ppw(G)$, is the minimum width over all proper-path decompositions of $G$.
If $\mathcal{X}$ satisfies conditions $1$-$4$ above, $\mathcal{X}$ is called a path-decomposition, which coincides with the definition of pathwidth given in Section~\ref{sec:intro}.

From the definitions, for any graph $G$ it clearly holds that
\begin{equation}
\ppw(G) \geq \pw(G) \geq \tw(G).
\end{equation}

Let us show that any element of $\Rcal_{n,k}$ has proper-pathwidth at most $k$.
Let $(\sigma, f, N)$ be constructible such that $G(\sigma,f,N) \in \Rcal_{n,k}$ and let $A_i$ be defined as in Subsection~\ref{sub:const}.
We define for every $i \in \{k+1,\ldots, n\}$ the bag $X_i = \{\vv{j} \mid j \in A_i \cup \{i\}\}$.
The sequence $\mathcal{X} = \{X_{k+1},X_{k+2}, \ldots, X_n\}$ satisfies the five conditions of the above definition, and
 for every $i \in \{k+1,\ldots, n\}$, $|X_i| = k+1$.
It follows that $G(\sigma,f,N)$ has proper-pathwidth at most $k$,
so it also has treewidth at most $k$, and therefore $G(\sigma,f,N)$ is a partial $k$-tree by Lemma~\ref{lemma:tw}.

\section{Counting the number of elements}
\label{sec:count}

In this section we analyze our construction and give a lower bound on $R_{n,k}$. We first start in
Subsection~\ref{sub:numb} by counting the number of constructible triples $(\sigma, f, N)$ generated by the algorithm, and in Subsection~\ref{sub:bound} we provide an upper bound on the number of duplicates. Finally, in Subsection~\ref{sub:s} we argue about the best choice for
 the parameter $s$ defined in the construction.

\subsection{Number of constructible triples  $(\sigma, f, N)$}
\label{sub:numb}

We proceed to count the number of constructible triples $(\sigma, f, N)$
created by  the construction given in Subsection~\ref{sub:const}.
As $\sigma$ is a permutation of $\{1,\ldots,n\}$ with the constraint that $\sigma(1) = 1$,
there are $(n-1)!$ distinct possibilities for the choice of  $\sigma$.
The function $f$ can take more than one value only for $k+2 \leq i \leq n$ and  $i \not\equiv {k+2 \pmod s}$.
This represents $n - (k+1) - \lceil\frac{n - (k+1)}{s}\rceil$ cases.
In each of these cases, there are at least $\frac{k-1}{2}$ distinct possible values for $f(i)$.
Thus, we have at least $(\frac{k-1}{2})^{(n - (k+1) - \lceil\frac{n - (k+1)}{s}\rceil)}$ distinct possibilities
for the choice of $f$.
For every $i \in \{2, \ldots, k+1\}$, $N(i)$ can be chosen as any subset of $i-1$ vertices containing the fixed vertex $\vv{1}$.
This yields $\prod_{i=2}^{k+1} 2^{i-2} = 2^{\frac{k(k-1)}{2}}$ ways to define $N$ over $\{2, \ldots, k+1\}$.
For $i \geq k+2$, $N(i)$ can be chosen as any subset of size at least $\frac{k+1}{2}$ of a set of $k$ elements with one element that is imposed.
This results in $\sum_{i=\lceil\frac{k+1}{2}\rceil}^k {k-1 \choose i-1} \geq 2^{k-2}$ possible choices for $N(i)$.
Thus, we have $2^{\frac{k(k+1)}{2}} \cdot 2^{(n-(k+1))(k-2)}$ distinct possibilities to construct $N$.

By combining everything, we obtain at least
\begin{equation}
(n-1)! \cdot \left(\frac{k-1}{2}\right)^{n - (k+1) - \lceil\frac{n - (k+1)}{s}\rceil} \cdot 2^{\frac{k(k-1)}{2}} \cdot 2^{(n-(k+1))(k-2)}
\end{equation}
distinct possible constructible triples $(\sigma, f, N)$.

\subsection{Bounding the number of duplicates}
\label{sub:bound}

Let $H$ be an element of $\Rcal_{n,k}$.  Our objective is to obtain an
upper bound on the number of constructible triples $(\sigma,f,N)$ such
that $H = G(\sigma, f, N)$.

Given $H$, we start by reconstructing $\sigma$.  Firstly, we
know by construction that $\sigma(1) = 1$.  Secondly, we know that $f(k+2) = 1$ and
so, for every $i > k+1$, $1 \not \in A_i$, implying that
$1 \not \in N(i)$.  It follows that the only neighbors of $\vv{1}$ are
$\vv{i}$ with $ 1 < i \leq k+1$.
So the set of images by $\sigma$ of $\{2, \ldots, k+1\}$ is uniquely
determined.  Then we guess the function $\sigma$ over this set
$\{2, \ldots, k+1\}$.  We have $k!$ possible such guesses for $\sigma$.

\begin{figure}[h]
\centering
\begin{tikzpicture}
  \foreach \i in {1,...,5}
  {
    \node[place,label=left:$v_{\sigma(i_{\i})}$] at (0,5-\i) {};
  }

\draw (-1,4.5) -- (-1.5,4.5) -- (-1.5,-0.5) -- (-1,-0.5);
\capt{-3.5,2}{$A_{i-1}$};

  \foreach \i in {0,...,10}
  {
    \node[place] at (2+0.5*\i,1) {};
  }
\draw (1.75,1) -- (1.75,0.5) -- (3.75,.5) -- (3.75,1);

  \foreach \i in {0,...,3}
  {
    \draw (0,4) -- (2+0.5*\i,1);
  }
  \foreach \i in {4,...,10}
  {
    \draw[dotted] (0,4) --  (2+0.5*\i,1) {};
  }

\anch{0,4}

\capt{1.7,0}{block of $s$ vertices};
\end{tikzpicture}
  \caption{\ The current anchor $\vv{i_1}$ is connected to all the $s$ vertices of the current block but will not be connected to any of the remaining non-introduced vertices.
}
  \label{fig:dupli}
\end{figure}
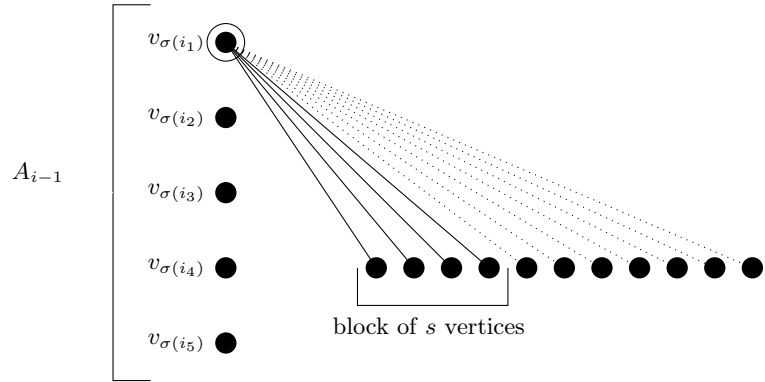

Thirdly, assume that we have correctly guessed $\sigma$ on
$\dotsss{k+1+ps}$ for some non-negative integer $p$ with $k+1+ps < n$.  Then
$a_{k+1+ps+1}$ is the smallest active vertex that is adjacent to at
least one element that is still not introduced after step $k+1+ps$.
Then the neighbors of $a_{k+1+ps+1}$ over the elements that are
not introduced yet after step $k+1+ps$ are the elements whose indices are
between $k+1+ps+1$ and $k+1+(p+1)s$, and these vertices constitute the next block of the construction; see Fig.~\ref{fig:dupli} for an illustration.  As before, the set of images by
$\sigma$ of $\{k+1+ps+1, \ldots, k+1+(p+1)s\}$ is uniquely determined,
and we guess $\sigma$ over this set.  We have at most $s!$ possible such
guesses.
Fourthly, if $ k+1+(p+1)s > n$ (that is, for the last block, which may have size smaller than $s$), we have $t!$ possible guesses with
$t = n -(k+1) - s\lfloor\frac{n-(k+1)}{s}\rfloor$.

We know that the first, the second, and the fourth cases can occur only once in the construction, and
the third case can occur at most $\lfloor\frac{n-(k+1)}{s}\rfloor$
times.  Therefore, an upper bound on the number of distinct possible guesses
for $\sigma$ is
$k!\cdot (s!)^{\lfloor\frac{n-(k+1)}{s}\rfloor} \cdot t!$, where
$t = n -(k+1) - s\lfloor\frac{n-(k+1)}{s}\rfloor$ .

Let us now fix $\sigma$.  Then the function $N$ is uniquely
determined.  Indeed, for every $i \in \dotsss{n}$, $N(i)$ corresponds
to the neighbors of $\vv{i}$ to the left.  It remains to bound the
number of possible functions $f$.  In order to do this, we define for
every $i > 1$,
$D_i = \{j \in N(i) \mid \forall j' > i, \{\vv{j}, \vv{j'}\} \not \in
E(H)\}$.  Then, for every $i \geq k+2$, by definition of $f(i)$,
$f(i) \in D_{i-1}$.  Moreover, for $i,j > k+1$ with $ i \not = j$, it
holds that, by definition of $D_i$ and $D_j$, $D_i \cap D_j = \es$.
Indeed, assume w.l.o.g. that $i < j$, and suppose for contradiction that there exists $a \in D_i \cap D_j$.
As $a \in D_j$, it holds that $a \in N(j)$, but as $a \in D_i$,  for every $j' > i$, $a \not \in N(j')$, hence $a \not \in N(j)$, a contradiction.

We obtain that the number of distinct functions $f$ is bounded by
$\prod_{i=k+1}^{n}|D_i|$.  As $D_i \cap D_j = \es$ for every $i,j \geq k+1$
with $i \not = j$ and $D_i \subseteq \dotsss{n}$ for every $i \geq k+1$,
we have that $\sum_{i=k+1}^{n} |D_i| \leq n$.
Let $I = \{i \in \{k+1, \ldots, n\} \mid |D_i| \geq 2\}$, and note that $|I| \leq k$.
By the previous discussion, it holds that $\sum_{i \in I} |D_i| \leq 2k$.
So it follows that, by using Cauchy-Schwarz inequality,
\begin{equation}
\prod_{i=k+1}^n |D_i|\ =\ \prod_{i \in I} |D_i|\ \leq\ \left(\frac{\sum_{i \in I}|D_i|}{k}\right)^k \ \leq\ \left(\frac{2k}{k}\right)^k\ =\ 2^k.
\end{equation}

To conclude, the number of constructible triples that can give rise to
$H$ is at most
$2^k \cdot (s!)^{\lfloor\frac{n-(k+1)}{s}\rfloor} \cdot t!$ where
$t = n -(k+1) - s\lfloor\frac{n-(k+1)}{s}\rfloor$.  Thus, we obtain
that

\begin{equation}
  \label{eq:rnk}
  R_{n,k}\ \geq\ \frac{(n-1)! \cdot \left(\frac{k-1}{2}\right)^{n - (k+1) - \lceil\frac{n - (k+1)}{s}\rceil} \cdot 2^{\frac{k(k-1)}{2}} \cdot 2^{(n-(k+1))(k-2)}}{2^k \cdot k!\cdot (s!)^{\lfloor\frac{n-(k+1)}{s}\rfloor} \cdot (n -(k+1) - s  \lfloor\frac{n-(k+1)}{s}\rfloor)!}.
\end{equation}

For better readability, we bound separately each of the terms of Equation~(\ref{eq:rnk}):

\begin{itemize}
\item[$\bullet$] $(n-1)! \geq \frac{1}{n}(\frac{n}{e})^n$.
\item[$\bullet$] $\frac{1}{n} \geq 2^{-n}$.
\item[$\bullet$]
  $({k-1})^{(n - (k+1) - \lceil\frac{n - (k+1)}{s}\rceil)} \geq
  2^{-n}k^{(n-\frac{n}{s} - k -2)}$, where we have assumed that $k \geq 2$, in which case $k-1 > \frac{k}{2}$; if $k=1$, we already know that $T_{n,1} \sim \sqrt{e} \cdot n^{n-2}$~\cite{Forest90}.
\item[$\bullet$] $2^{(n - (k+1) - \lceil\frac{n - (k+1)}{s}\rceil)} \leq 2^n$.
\item[$\bullet$]
  $2^{\frac{k(k-1)}{2}} \cdot 2^{(n-(k+1))(k-2)} \geq
  2^{kn-\frac{k(k+3)}{2}} \cdot 2^{-2n}$.
\item[$\bullet$] $2^k \leq 2^n$.
\item[$\bullet$] $k! \leq k^k$.
\item[$\bullet$]
  $(s!)^{\lfloor\frac{n-(k+1)}{s}\rfloor} \cdot (n -(k+1) -
  s\lfloor\frac{n-(k+1)}{s}\rfloor)! \leq s^n$.
\end{itemize}

By applying these considerations into Equation~(\ref{eq:rnk}), we can
simplify it to
\begin{equation}
  \label{eq:srnk}
  R_{n,k}\ \geq\ \left(\frac{1}{64e} \cdot \frac{k \cdot 2^{k} \cdot n }{ k^{\frac{1}{s}} \cdot s}\right)^n\cdot 2^{-\frac{k(k+3)}{2}} \cdot k^{-2k -2}.
\end{equation}


\subsection{Choice of the parameter $s$}
\label{sub:s}

Let us now discuss how to choose the size $s$ of the blocks in the construction.
In order to obtain the largest possible lower bound for $R_{n,k}$, we would like to choose the  value of $s$ that minimizes the denominator $k^{\frac{1}{s}} \cdot s$ in Equation~(\ref{eq:srnk}).
To be as general as possible,  assume that $s$ is a function $s(n,k)$ that may depend on $n$ and $k$, and  we define $t(n,k) := \frac{s(n,k)}{\log k}$.
With this definition, it follows that
\begin{equation}
\log \left(k^{\frac{1}{s(n,k)}} \cdot s(n,k)\right)\ =\ 	\frac{\log k}{s(n,k)}+ \log s(n,k)\ =\ \frac{1}{t(n,k)} + \log t(n,k) + \log \log k.
\end{equation}

It is elementary that the minimum of $\frac{1}{t(n,k)} + \log t(n,k)$ is achieved for $t(n,k) = 1$.
Thus, we obtain that $s(n,k)= \log k$ is the function that
maximizes the lower bound given by Equation~(\ref{eq:srnk}).
Therefore, we obtain that
\begin{equation}
  R_{n,k}\ \geq\ \left(\frac{1}{128e} \cdot \frac{k \cdot 2^{k} \cdot n }{ \log k}\right)^n\cdot 2^{-\frac{k(k+3)}{2}} \cdot k^{-2k -2},
\end{equation}
concluding the proof of Theorem~\ref{thm:main}, where we assume that $k \geq 2$.

\section{Concluding remarks and further research}
\label{sec:further}

Comparing Equation~(\ref{eq:sub}) and Equation~(\ref{eq:ilbi}), there is still a gap of
$(128e \cdot \log k)^n$ in the dominant term of $T_{n,k}$, and closing this gap remains a challenging open problem. The factor $(\log k)^n$ appears because, in our construction, when a new block starts, that is, every $s = \log k$ introduced vertices, we force the  frozen vertex to be the previous anchor. Therefore, this factor is somehow artificial, and we believe that it could be improved.


 One could also focus on the term of $T_{n,k}$ that depends only on $k$, namely $2^{-\frac{k(k+3)}{2}}\cdot k^{-2k -2}$ for the lower bound and $2^{ - \frac{k(k+1)}{2}} \cdot k^{-k}$ for the upper bound. 
In our lower bound, we think that the constant $3$ in the term $\frac{k(k+3)}{2}$ may be reduced to $1$,  as its existence is related  to the fact that, in the construction, we force $\sigma(1) = 1$, and therefore the neighborhood of the first $k+1$ vertices, except for the first one,  is forced to contain vertex $1$.

We believe that there exist an absolute constant $c > 0$ and a function $f(k)$, with $k^{-2k-2} \leq f(k) \leq k^{-k}$ for every $k>0$, such that for every $0 < k \leq n$,
\begin{equation}
T_{n,k}\ \geq\ (c\cdot k\cdot 2^k \cdot n)^n\cdot  2^{-\frac{k(k+1)}{2}} \cdot f(k).
\end{equation}

One way to improve the upper bound would be to show that {\sl every} partial $k$-tree with $n$ vertices and $m$ edges can be extended to at least a large number $\alpha(n,m)$ of $k$-trees, and then use \emph{double counting}. This is the approach taken in~\cite{OsthusPT03} for bounding the number of planar graphs, but  so far we have not been able to obtain a significant improvement using this technique.

%

\medskip

Our results find algorithmic applications, specially in the area of Parameterized Complexity. When designing a parameterized algorithm, usually a crucial step is to solve the problem at hand restricted to graphs decomposable along small separators by performing dynamic programming (see~\cite{JansenLS14} for a recent example).
For instance, precise bounds on $T_{n,k}$ are useful when dealing with the
\textsc{Treewidth-$k$ Vertex Deletion} problem, which has recently attracted significant attention in the area~\cite{LPRRSS16,GHOORRVS13,FominLMS12}. In this problem, given a graph $G$ and a fixed integer $k>0$, the objective is to remove as few vertices from $G$ as possible in order to obtain a graph of treewidth at most $k$. When solving \textsc{Treewidth-$k$ Vertex Deletion} by dynamic programming, the natural approach is to enumerate, for any partial solution at a given separator of the decomposition,  all possible graphs of treewidth at most $k$ that are ``rooted'' at the separator. In this setting, the value of $T_{n,k}$, as well as an explicit construction to generate such graphs, may be crucial in order to speed-up the running time of the  algorithm.

%
%

As mentioned before, our results also apply to other relevant graph parameters such as pathwidth and proper-pathwidth. For both parameters, beside improving the lower bound given by our construction, it may be also possible to improve the upper bound given by Equation~(\ref{eq:sub}). For proper-pathwidth, a modest such improvement can be obtained by improving the upper bound given by Theorem~\ref{thm:ktree}. Indeed, it easily follows from the definition of proper-pathwidth  that the edge-maximal graphs of proper-pathwidth $k$, which we call \emph{proper linear $k$-trees}, can be constructed starting from a $(k+1)$-clique and iteratively adding a vertex
$v_i$ connected to a clique $K_{v_i}$ of size $k$, with the constraints that $v_{i-1} \in K_{v_i}$ and $K_{v_i} \setminus \{v_{i-1}\} \subseteq K_{v_{i-1}}$. From this observation, and taking into account that the order of the first $k$ vertices is not relevant and that there are $2k$ initial cliques giving rise to the same graph, it follows that the number of $n$-vertex labeled proper linear $k$-trees is equal to
\begin{equation}\label{eq:number-proper-linear-k-trees}
n! \cdot k^{n-k-1} \cdot \frac{1}{2k \cdot k!}.
\end{equation}

From Equation~(\ref{eq:number-proper-linear-k-trees}) and using that an $n$-vertex labeled proper linear $k$-tree has $kn - \frac{k(k+1)}{2}$ edges, basic calculations yield that the dominant term of the number of $n$-vertex labeled graphs of proper-pathwidth at most $k$ is at most $\left(\frac{k \cdot 2^{k} \cdot n}{c} \right)^n$ for some absolute constant $c \geq 1.88$.

Finally, it would be interesting to count the graphs of bounded $X$\emph{width}, for other $X$ different than ``tree'', ``path'', or ``proper-path''. For instance, branchwidth seems to be a good candidate, as it is known that, if we denote by $\bw(G)$ the branchwidth of a graph $G$ and $|E(G)| \geq 3$, then $\bw(G) \leq \tw(G) +1 \leq \frac{3}{2}\bw(G)$~\cite{RoSe91}. Other relevant graph parameters are cliquewidth, rankwidth, tree-cutwidth, or booleanwidth. For any of these parameters, a first natural step would be to find a ``canonical'' way to build such graphs, as it is the case of partial $k$-trees.


\medskip

\noindent{\small\textbf{Acknowledgement}. We would like to thank Dimitrios M. Thilikos for pointing us to the notion of proper-pathwidth.}

\vspace{-.25cm}

{\small
\bibliographystyle{abbrv}
\bibliography{biblio-tw}
}

\end{document}

%% file: newcommand.tex
\newcommand{\Rcal}{\mathcal{R}}

\newcommand{\sm}{\setminus}

\newcommand{\tw}{{\mathbf{tw}}}
\newcommand{\bw}{{\mathbf{bw}}}
\newcommand{\pw}{{\mathbf{pw}}}
\newcommand{\ppw}{{\mathbf{ppw}}}

\newcommand{\dotsss}[1]{\{1,\dots, #1\}}

\newcommand{\es}{\emptyset}

\newcommand{\vv}[1]{v_{\sigma(#1)}}

\newcommand{\capt}[2]{\draw[white] (#1) -- node[above,sloped,black] {#2} +(2,0)}

\tikzstyle{place}=[
  circle,
  fill=black,
  thick,
  inner sep=0.1cm,
  minimum size=1mm
]

\newcommand{\forg}[1]{\draw (0.3+#1+0.3) -- (-0.3+#1-0.3);\draw (0.3+#1-0.3) -- (-0.3+#1+0.3);}
\newcommand{\anch}[1]{\draw (#1) circle (2.5mm);}